\documentstyle{article}
\begin{document}
\newtheorem{thm}{Theorem}[section]
\newtheorem{lemma}[thm]{Lemma}
\newtheorem{defin}[thm]{Definition}
\newtheorem{rem}[thm]{Remark}
\newtheorem{cor}[thm]{Corollary}
\newtheorem{prop}[thm]{Proposition}

\def\fn{ \baselineskip = 0pt\vbox{\hbox{\hspace*{3pt}\tiny $\circ$}\hbox{$f$}}\baselineskip = 12pt\!}
\def\gn{ \baselineskip = 0pt\vbox{\hbox{\hspace*{2pt}\tiny $\circ$}\hbox{$g$}}\baselineskip = 12pt\!}
\def\mn{ \baselineskip = 0pt\vbox{\hbox{\hspace*{3pt}\tiny $\circ$}\hbox{$m$}}\baselineskip = 12pt\!}
\def\lap{\bigtriangleup}
\def\be{\begin{equation}}
\def\ee{\end{equation}}
\def\bea{\begin{eqnarray}}
\def\eea{\end{eqnarray}}
\def\beas{\begin{eqnarray*}}
\def\eeas{\end{eqnarray*}}
\def\dt{\partial_t}
\def\dx{\partial_x}
\def\dv{ \partial_v }
\def\dtau{\partial_\tau}
\def\dxi{\partial_\xi}
\def\deta{ \partial_\eta }
\def\tf{\tilde f}

\def\R{{\rm I\kern-.1567em R}}
\def\N{{\rm I\kern-.1567em N}}
\def\Z{{\sf Z\kern-.3567em Z}}
\def\ekin{E_{\rm k}}
\def\epot{E_{\rm p}}
\def\supp{\mbox{\rm supp}}
\def\div{\mbox{\rm div}}
\def\n#1{\vert #1 \vert}
\def\nn#1{{\Vert #1 \Vert}}
\def\prf{\noindent{\bf Proof.\ }}
\def\prfe{\hspace*{\fill} $\Box$ \smallskip \noindent}
\def\eqn#1{(\ref{#1})}

\def\un {{\rm 1{\kern-.26em }I}}
\def\wrt{with\ respect\ to\ }

\def\today{\space\ifcase\month\or January\or February \or March \or April
\or May\or June\or July\or August\or September\or October\or November\or
December\fi\space \number\day ,\space\number\year}

\newcommand{\square}{$\Box$}
\def\ix#1{\int{#1}\; dx}\def\ixd#1{\int_{\R^d}{#1}\; dx} 
\def\ix2#1{\int_{\R^2}{#1}\; dx} 
\def\ixv#1{\int\!\!\!\!\!\int{#1}\; dxdv}
\def\ixe#1{\int\!\!\!\!\!\int{#1}\; dxd\eta} 
\def\ixed#1{\int\!\!\!\!\!\int_{\R^d\times\R^d}{#1}\; dxd\eta} 
\def\ixe2#1{\int\!\!\!\!\!\int_{\R^2\times\R^2}{#1}\; dxd\eta} 
\def\ixi#1{\int{#1}\; d\xi}
\def\ixie#1{\int\!\!\!\!\!\int{#1}\; d\xi d\eta}
\def\rn{{\R^d}}

\title{Time-dependent rescalings and Lyapunov functionals for
the Vlasov-Poisson and Euler-Poisson systems, and for related models 
of kinetic equations, fluid dynamics and quantum physics}
\author{Jean Dolbeault\\
Ceremade, U.M.R. C.N.R.S. no. 7534,\\
Universit\'e Paris IX-Dauphine\\
Place du Mar\'echal de Lattre de Tassigny\\
75775 Paris C\'edex 16, France\\
and\\
Gerhard Rein\\ 
Mathematisches Institut, Universit\"at M\"unchen\\
Theresienstr. 39, D-80333 M\"unchen, Germany}
\date{\today}
\maketitle

\begin{abstract}
We investigate rescaling transformations for the Vlasov-Poisson and
Euler-Poisson systems and derive in the plasma physics case Lyapunov 
functionals which can be used to analyze dispersion effects. The method is
also used for studying the long time behaviour of the solutions and can be
applied to other models in kinetic theory (2-dimensional symmetric
Vlasov-Poisson system with an external magnetic field), in fluid dynamics
(Euler system for gases) and in quantum physics (Schr\"odinger-Poisson system,
nonlinear Schr\"odinger equation). 
\end{abstract}\bigskip

\bigskip

\noindent {\bf 1991 MR Subject Classification.} 
35B40,
82C40.
\vfill\eject

\section{Introduction}
\setcounter{equation}{0}
Consider the Vlasov-Poisson system (VP)
\[ 
\dt f + v\cdot \dx f - \dx U \cdot \dv f =0, \]
\[ \lap U = \varepsilon \rho,\
\rho (t,x) = \int_{\R^d} f(t,x,v)\,dv 
\]
and the pressureless Euler-Poisson system (EP)
\[ 
\dt \rho + \div (\rho \,u) = 0, 
\]
\[ 
\dt u + (u \cdot \dx) u = - \dx U, 
\]
\[ 
\lap U = \varepsilon \rho . 
\]
Here $t \geq 0$, $x$, $v$, $u=u(t,x) \in \R^d$, $d\geq 1$ is the dimension of 
the physical space, $\varepsilon = + 1$ corresponds to the stellar dynamics 
and $\varepsilon = -1$ to the plasma physics case.
Throughout this paper, we shall assume that $f$ is a nonnegative function
in $L^\infty(\R^+,\, L^1(\R^d\times\R^d))$. 
Formally, we have the following relation between these two systems:
a pair $(\rho,u)$ is a solution of (EP) if and only if
\[ 
f(t,x,v)= \rho (t,x) \delta\Bigl(v - u(t,x)\Bigr) 
\]
is a solution of (VP) where $\delta$ denotes the Dirac delta distribution.
In this situation $u$ can be recovered from $f$ via the identity
\be\label{VPEP}
\rho (t,x)\; u(t,x) = \int_{\R^d} v\, f(t,x,v)\, dv .
\ee
In this sense (EP) is a special case of (VP), and we will see later that the
asymptotic behaviour of (VP) for large times is connected with a
special solution of (EP). On a rigorous level the relation of (VP) with (EP)
is investigated in \cite{Di}.

\medskip
Throughout this paper, we assume for simplicity that the solutions of (VP) are
of class $C^1$ with compact support with respect to
$x$ and $v$, which allows us to perform any integration
by parts without further justifications (except maybe in dimension $2$). The
results then pass to less smooth classes of solutions, assuming for instance
that $f$ belongs to
$C^0(\R^+,L^1(\R^d\times\R^d))$ (see for instance \cite{P1} or \cite{Sc}) and is
a global in time solution to the Cauchy problem corresponding to an initial
data $f_0$ satisfying for instance:
\begin{description}
\item{($d=2$)}
$f_0\in L^1\cap L^\infty(\R^2\times\R^2)$ 
is such that for some $\epsilon>0$ the quantity
\[ 
\int\!\!\!\!\!\int_{\R^2\times\R^2}f_0(x,v)\; (\vert x\vert^{2+\epsilon}+\vert
v\vert^{2+\epsilon}+\vert U_0(x)\vert)\; dxdv 
\]
(with $U_0(x)=-\frac 12 \log\vert x\vert *\int f_0\, dv$) is bounded 
(see \cite{D1}).
\item{($d=3$)}
$f_0\in L^1\cap L^\infty(\R^3\times\R^3)$
is such that for some $\epsilon>0$ and $p>3$ the quantity
\[ 
\int\!\!\!\!\!\int_{\R^3\times\R^3}f_0(x,v)(\vert v\vert^{2+\epsilon}
+\vert x\vert^p)\; dxdv 
\]
is bounded (see \cite{LP,Pf,R,Sc} and
\cite{C3,CP} for the propagation of moments).
\end{description}
For weak solutions obtained as a limit of an approximating sequence (for
instance, if we assume no moments higher than $2$), the equalities have to be
replaced by inequalities.
\medskip 

For the Euler-Poisson system, we shall consider only $C^1$ solutions.
The results presented in this paper have to be understood as either
a general method on how to obtain dispersion effects without 
taking care of the existence or the regularity of the solutions, or
as a method to derive {\sl a priori\/} estimates for less regular solutions
(by passing to the limit with smooth approximating solutions).

Our paper is organized as follows. In Section 2, we introduce linear 
scalings and explain why they give rise to singular self-similar problems.
How to remedy this pathology with time-dependent scalings is explained 
in Section~3. In the one-dimensional case, the information on the solution 
is sufficient to provide the convergence of the rescaled solution to the 
asymptotic measure. Section 4 is concerned with the Lyapunov functionals 
and constitutes the heart of this paper: the energy of the rescaled system 
turns out to be a Lyapunov functional for the initial problem. A more 
straightforward (than the full time-dependent scaling method) approach 
to the Lyapunov functionals is also given. In Section 5, we use the Lyapunov
functionals to describe the asymptotic behaviour (dispersion rate) of the solutions in the plasma physics case. 

Rescalings for the study of large time behaviour have been widely used in
various fields of applied mathematics but appear to be rather new in the
context of kinetic equations: in that direction we may mention the studies
made by J.~R.~Burgan, M.~R.~Feix, E.~Fijalkow, and A.~Munier (see
\cite{BFFM}) and J.~Batt, M.~Kunze, G.~Rein in \cite{BKR}. Our main point is to
make the link between rescalings preserving the $L^1$-norm and Lyapunov
functionals (or pseudo-conformal laws) and to explain on various examples of
conservative systems why it actually provides a general method for the study of
large time asymptotics.

While Sections 1--5 are exclusively devoted to the Vlasov-Poisson and 
Euler-Poisson systems, Sections 6--8 are concerned with other problems 
of kinetic theory, fluid mechanics and quantum physics. The relation 
between these various domains has been noticed for a long time (see for 
instance \cite{P2}), but it is surprising that the 
estimates given in \cite{IZL} have been adapted to 
kinetic models only recently. 
Here we proceed in the reverse historical order, from kinetic 
equations to fluids and quantum physics, and this approach actually 
seems to be very powerful.

To conclude with the introduction, it is worth mentioning that many of the 
estimates we are giving in this paper were already at least partially
known. The point is that we present a systematic and elementary method
which takes the nonlinearity of the model very well into account (this was
not necessarily the case in the preceding papers) and gives rise to a
more precise form of the Lyapunov functionals (in the sense that these
Lyapunov functionals also include second-moments in the $x$-variable) which
are natural for the problems we consider.

\section{Linear scalings}
\setcounter{equation}{0}

Let $f = f(t,x,v)$ be a solution of (VP). Then for any $\lambda,\mu >0$
\[ 
f_{\lambda,\mu} (t,x,v) =
\lambda^{2-d} \mu^d f(\lambda t,\mu x,\lambda^{-1} \mu v) 
\]
is again a solution of (VP),
\[ 
\dt f_{\lambda,\mu} + v\cdot \dx f_{\lambda,\mu} - \dx U_{\lambda,\mu} 
\cdot \dv f_{\lambda,\mu} =0, 
\]
\[ 
\lap U_{\lambda,\mu} = \varepsilon \rho_{\lambda,\mu},\quad
\rho_{\lambda,\mu} (t,x) = \int_{\R^d} f_{\lambda,\mu}(t,x,v)\,dv 
\]
with
\beas
\rho_{\lambda,\mu}(t,x)
&=&
\lambda^2 \rho(\lambda t,\mu x),\\
U_{\lambda,\mu} (t,x)
&=&
\lambda^2 \mu^{-2} U(\lambda t,\mu x),\\
\dx U_{\lambda,\mu} (t,x)
&=&
\lambda^2 \mu^{-1} \dx U(\lambda t,\mu x),
\eeas
as can be checked by direct computation using for instance the
following integral representation of $\partial_xU$:
\[ 
\partial_xU(t,x)=-\frac
x{\vert S^{d-1}\vert\, \vert x\vert^d}*\int_{\R^d} f(t,x,v)\; dv . 
\]
Here $S^{d-1}\subset \R^d$ is the unit sphere of dimension $d-1$.
Similarly, if $(\rho,u) = (\rho,u)(t,x)$ is a solution of (EP) then
for any $\lambda,\mu >0$
\[ 
\rho_{\lambda,\mu} (t,x) =
\lambda^2 \rho(\lambda t,\mu x),\
u_{\lambda,\mu} (t,x) = \lambda \mu^{-1} u(\lambda t,\mu x) 
\]
is again a solution of (EP), and the potential is transformed as for (VP).
This also follows from Relation (\ref{VPEP}) between (VP) and (EP).

If we require that the $L^1$-norm of $\rho (t)$,
which is a conserved quantity for (EP) as well as for (VP), is preserved 
by the scaling, $\lambda$ and $\mu$ must satisfy
\[ 
\lambda^2\mu^{-d}=1 
\]
and the rescaled distribution function is $f_{\lambda,\lambda^{2/d}}$.
A standard way of studying the asymptotic behaviour of $f$ would then be to 
consider a self-similar solution, {\sl i.e.\/} a solution which
satisfies
\[ 
f(t,x,v)=f_{\lambda,\lambda^{2/d}}(t,x,v)=\lambda^{4-d}f(\lambda t,
\lambda ^{2/d}x, \lambda ^{2/d-1}v) 
\]
for any $\lambda>0$. This solution would then be given by its self-similar
profile $\tilde f (\xi,\eta) = f(1,\xi,\eta)$ (choose $\lambda$ to be
$\frac 1t$). Then 
\be\label{sssln}
f(t,x,v)=f_{t^{-1}, t^{-2/d}}(t,x,v)
=t^{d-4}\tilde f (t^{-\frac 2d}x,t^{1-\frac 2d}v )
\ee
is a solution of (VP) if, at least formally, $\tilde f$ is a solution of
\[ 
(d-4)\tilde f + \eta\cdot \partial_\xi\tilde f 
- \frac 2d \xi\cdot \partial_\xi\tilde f
+(1-\frac 2d) \eta\cdot \partial_\eta\tilde f
- \partial_\xi\tilde U \cdot \partial_\eta\tilde f =0, 
\]
\[ 
\lap_\xi\tilde U = \varepsilon \tilde \rho,\quad
\tilde \rho (\xi) = \int_{\R^d}\tilde f(\xi,\eta)\,d\eta 
\]
in the new variables $\xi=t^{-\frac 2d}x$, $\eta =t^{1-\frac 2d}v$.
However, it is clear that as $t\to 0+$, $g(t,x,v)$ does not converge to 
a well defined measure for which one might establish an 
existence result, except for $d=1$. 
This difficulty is completely removed by considering 
general, non-singular, time-dependent scalings.

\section{Time-dependent scalings}
\setcounter{equation}{0}

Consider the following transformation of variables in (VP),
where the positive functions $A(t),\ R(t),\ G(t)$ will be determined later:
\[ 
dt = A^2(t) d\tau,\quad x = R(t) \xi. 
\]
Thus, assuming that $t\mapsto x(t)$ and $\tau\mapsto \xi(\tau)$ 
satisfy $\frac{dx}{dt}=v$ and $\frac{d\xi}{d\tau}=\eta$ respectively,
the new velocity variable $\eta$ has to satisfy
\[ 
v = \frac{dx}{dt} = \dot R(t) \xi + R(t) \frac{d\xi}{d\tau} \frac{d\tau}{dt}
= \dot R(t) \xi + \frac{R(t)}{A^2(t)} \eta. 
\]
Let $F$ be the rescaled distribution function:
\[ 
f(t,x,v) = G(t) F(\tau, \xi,\eta). 
\]
The aim is to choose this transformation in such a way that the transformed
Vlasov equation is still a transport equation on phase space and
contains a given, external force and a friction term.
The inverse transformation is
\[ 
d\tau = A^{-2}(t) dt,\ \xi = R^{-1}(t) x,\
\eta = \frac{A^2(t)}{R(t)} \biggl(v - \frac{\dot R (t)}{R (t)} x\biggr) \; . 
\]
Here $\dot {\phantom x}$ always denotes derivative \wrt $t$.
If $\nu$ and $W$ are defined as the rescaled spatial density
and the rescaled potential respectively, then
\[ 
\nu (\tau,\xi) = \int_{\R^d}F(\tau,\xi,\eta)\; d\eta = \frac{A^{2d}}{R^d G}
\rho (t,x), 
\]
\[ 
W(\tau,\xi) = \frac{A^{2d}}{R^{d+2} G} U (t,x),\quad
\dxi W(\tau,\xi) = \frac{A^{2d}}{R^{d+1} G} \dx U (t,x), 
\]
and the Vlasov equation transforms into
\beas
\dtau F + \eta \cdot \dxi F 
&+&
2 A^2 \biggl(\frac{\dot A}{A}
- \frac{\dot R}{R}\biggr) \eta \cdot \deta F \\
&-&
\ddot R \frac{A^4}{R} \xi \cdot \deta F
- R^d G A^{4-2d} \dxi W \cdot \deta F
+ A^2 \frac{\dot G}{G} F = 0 .
\eeas
We want this to be a conservation law on $(\xi,\eta)$-space, so we require
\be
\frac{\dot A}{A} - \frac{\dot R}{R} = \frac{1}{2d} \frac{\dot G}{G}
\ee
which holds if and only if
\be \label{cons}
G = c\,\cdot \left(\frac{A}{R}\right)^{2d}
\ee
for some constant $c >0$; recall that $G$ should be positive.
For simplicity and without any loss of generality, we may choose $c =1$,
and the Vlasov equation becomes
\[ \dtau F + \eta \cdot \dxi F
+\div_\eta \left[ \left(\frac{1}{d} A^2 \frac{\dot G}{G} \eta -
\ddot R \frac{A^4}{R} \xi - R^d G A^{4-2d} \dxi W\right) F\right] = 0 . \]
Next we require that the external force in the above Vlasov equation
becomes time-independent and that there is no time-dependent factor
in front of the nonlinear term, {\sl i.e.}
\be\label{simp1}
\ddot R \frac{A^4}{R} = - \varepsilon c_0,
\ee
\be\label{simp2}
R^d G A^{4-2d} = 1 ,
\ee
where $c_0 >0$ is an arbitrary constant.
In view of (\ref{cons}) we get
\[ A = R^{d/4},\ G = R^{\frac{d-4}{2}d} \]
and $R$ has to solve
\be\label{greq} 
\ddot R + \varepsilon c_0 R^{1-d} = 0 .
\ee
At this point, we may state the following remark: 

\begin{rem}\label{ODE} 
Every solution of Equation (\ref{greq})
has the following properties:
\begin{description}
\item{(i)}
For any $\lambda>0$,
$t\mapsto R_\lambda(t)=c_0^{-\frac 1d} \lambda^{-\frac 2d} 
R(\lambda t)$
is a solution of 
\be\label{req}
\ddot R + \varepsilon R^{1-d}=0\; .
\ee
Without loss of generality we therefore assume that $c_0=1$ in what follows.
\item {(ii)}
With $R_0=R(0)$ and $\dot{R}_0=\dot{R}(0)$ we get, for $d=1$,
\[ 
R(t)=-\frac \varepsilon 2 t^2+\dot{R}_0 t+R_0\; . 
\]
If $d\geq 2$, it is easy to carry out one integration of Equation (\ref{req}):
\[
\frac 12\dot{R}^2(t)+\varepsilon\log R(t) =\frac 12\dot{R}_0^2+\varepsilon\log R_0
\ \ \mbox{for}\ d=2,
\]
\[
\frac 12\dot{R}^2(t)-\frac\varepsilon{d-2} R^{2-d}(t)
=\frac 12\dot{R}_0^2-\frac\varepsilon{d-2} R^{2-d}_0\ \ \mbox{for}\ d\geq 3 .
\]
In the plasma physics case, $R(t)$ cannot change sign and is
well defined for any $t\in\R$. Moreover,
\[
\log R(t) =\frac 12\dot{R}^2(t)-\frac 12\dot{R}_0^2+\log R_0
\geq-\frac 12\dot{R}_0^2+\log R_0
\ \ \mbox{for}\ d=2,
\]
\[
0 \leq R^{2-d}(t)=
\frac{d-2}2\biggl(\dot{R}_0^2-\dot{R}^2\biggr) + R^{2-d}_0
\leq\frac{d-2}2\dot{R}_0^2 + R^{2-d}_0 
\ \ \mbox{for}\ d\geq 3 .
\]
Together with Equation (\ref{req}) this proves that  there
exists a unique
$t_0\in\R$ such that $R(t_0)>0$ and $\dot{R}(t_0)=0$, and $R(t)>R_0$ for any
$t\neq t_0$, provided $\varepsilon = -1$.
\item{(iii)}
If $t\mapsto R(t)$ is a solution of Equation (\ref{req}) with $\varepsilon=-1$,
$t\mapsto R(t+a)$ is a solution too for any given $a\in\R$. Combining this
with the invariance through the rescaling $\lambda\mapsto R_\lambda(t)$ with
$R_\lambda(t)=\lambda^{-\frac 2d}R(\lambda t)$, we may always require 
$R_0=1$ and $\dot{R}_0=0$ without loss of generality as long as we are 
interested in the asymptotic behaviour of $f$ when $t\to +\infty$. Note that 
with this special choice for $R_0$ and $\dot{R}_0$, at $t=0$, $G(0)=A(0)=1$, 
and if we assume $\tau(0)=0$, then
\[ 
\xi(\tau=0,x)=x,\quad\eta(\tau=0,x,v)=v\quad\mbox{and}\quad 
f(t=0,x,v)=F(\tau=0,x,v)\; . 
\]
The time-dependent rescaling has the interesting property that it does
not introduce any singularity at $t=0$, 
and with $R_0=1$ and $\dot{R}_0=0$, the
initial data for $f$ and $F$ are the same.
\item{(iv)}
The singular self-similar solution (\ref{sssln})
corresponding to the linear scalings of Section 2 is---when it exists---the solution one expects to get in the limit case $R_0=0$. 
Formally, this solution also corresponds to the
limit of $R_\lambda(t)$ as $\lambda\to +\infty$.
\item{(v)} For $\varepsilon=-1$ and $d\geq 1$, $ \frac{\dot{R}}R\sim\frac 1t$ as
$t\to +\infty$ and 
\be\label{Requiv}
\matrix{&R(t)\sim t^2\qquad &\mbox{for}\ d=1\;,\cr
&R(t)\sim t\sqrt{\log t}\qquad &\mbox{for}\ d=2\;,\cr
&R(t)\sim t\qquad &\mbox{for}\ d\geq 3\; .\cr}
\ee
\end{description}
\medskip
\end{rem}

With $R$ solving Equation \eqn{req}, we obtain the following rescaled
Vlasov-Poisson system (RVP):
\[
\dtau F + \eta \cdot \dxi F
+\div_\eta
\left[ \left( \varepsilon \xi - \dxi W
+ \frac{d-4}{2} R^{\frac{d}{2} -1} \dot R \eta \right) F \right]
= 0,
\]
\[
\lap W = \varepsilon \nu (\tau,\xi) = \varepsilon \int_{\R^d} F (\tau,\xi,\eta)\,
d\eta .
\]
The relation between the old and the new variables is
\[ 
dt = R^{d/2} d\tau,\quad d\tau = R^{- d/2} dt, 
\]
\[ 
x = R\, \xi,\quad \xi = R^{-1} x, 
\]
\[ 
v= \dot R \xi + R^{1-\frac{d}{2}} \eta,\quad
\eta = R^{\frac{d}{2}-1} \biggl( v - \frac{\dot R}{R} x\biggr), 
\]
and the rescaled functions are given by
\[
F (\tau,\xi,\eta) =
R^{\frac{4-d}{2}d} f(t,x,v),
\]
\[
\nu (\tau,\xi) =
R^d \rho (t,x),
\]
\[
W (\tau,\xi) =
R^{d-2} U(t,x),\quad \dxi W (\tau,\xi) =
R^{d-1} \dx U(t,x).
\]
If we consider (EP) we find that with
\[ 
\eta(\tau,\xi) = R^{\frac{d}{2} -1}
\biggl( u(t,x) - \frac{\dot R}{R} x \biggr) 
\]
the rescaled Euler-Poisson system (REP) is
\[ \dtau \nu + \div (\nu \,\eta) = 0, \]
\[ \dtau \eta + (\eta \cdot \dxi)\,\eta =
\varepsilon \xi - \dxi W +
\frac{d-4}{2} R^{\frac{d}{2} -1} \dot R \eta, \]
\[ \lap W = \varepsilon \nu . \]
Note that this rescaling as the one for the Vlasov-Poisson system introduces a
harmonic force term $\varepsilon\xi$ and a friction term which is
proportional to the velocity.\medskip 

There exists a unique steady state with a given $L^1$-norm $M$ for which the
particles are at rest and uniformly distributed in the unit ball centered
at $0$, and the self-consistent force is exactly balanced by the
external force. Define
\[ 
F_\infty^M (\xi,\eta) = \nu_\infty^M (\xi) \delta (\eta) 
\]
where $\delta$ is the usual Dirac distribution. If $F_\infty^M$
and $(\nu_\infty^M,\eta_\infty^M=0)$ are the stationary
solutions of (RVP) and (REP)  respectively such that
\[ 
\Vert F_\infty^M\Vert_{L^1(\R^d\times\R^d)}
=\Vert\nu_\infty^M\Vert_{L^1(\R^d)}=M, 
\] 
then
\be\label{asfield}
\dxi W_\infty^M (\xi) =
\varepsilon 
\left\{
\begin{array}{ccl}
\xi &,& \n{\xi} \leq (M/\vert S^{d-1}\vert)^{1/d},\\
\xi/\n{\xi}^d &,& \n{\xi} >(M/\vert S^{d-1}\vert)^{1/d},
\end{array} \right.
\ee
and
\be\label{asdens}
\nu_\infty^M (\xi) = d\,\cdot\, \un_{B^d((M/\vert
S^{d-1}\vert)^{1/d})}.
\ee
Here $B^d(r)$ denotes the ball with radius $r$ centered at $0\in\R^d$, and
$\un_\omega$ denotes the characteristic function of the set $\omega$.
The inverse rescaling transformation takes this steady state
into
\be\label{assvp}
f_\infty^M (t,x,v) = \frac{d}{R(t)^d} \un_{B^d(R(t)( M/\vert
S^{d-1}\vert)^{1/d})}(x)\,
\delta \biggl(v - \frac{\dot R (t)}{R(t)} x \biggr),
\ee
and
\[ 
\rho_\infty^M (t,x) = \frac{d}{R(t)^d} \un_{B^d(R(t)( M/\vert
S^{d-1}\vert)^{1/d})}(x),\quad
u_\infty^M (t,x) = \frac{\dot R (t)}{R(t)} x . 
\]
It is easy to see that this defines a weak solution of (VP)
or (EP) respectively.

\medskip
In the plasma physics case we have $\dot R(t) >0$ for $t>0$ provided $\dot R
(0) \geq 0$ so that for $d \leq 3$ the particles are slowed down by a
friction force, and on physical grounds one would expect that the steady state
written above is a global attractor for (RVP) or (REP) respectively. In
\cite{BKR} this was carried out rigorously for the case
$d=1$. We will see in the next section that this is not true in general, at
least in dimension $d=3$ for (EP). However, the rescaling still provides
informations on the asymptotic behaviour of the original system for large
times: the energy for the rescaled system gives rise to a
Lyapunov functional for the original system by which dispersion
effects and the asymptotic behaviour can be analyzed.

\section{Lyapunov functionals}
\setcounter{equation}{0}

In this section, we investigate the behaviour of the total energy of
(RVP) and (REP). Let us consider first the case $d\geq 3$. 
The potential energy term is the same for both systems, namely
\[ \epot (\tau) = \int ( W(\tau,\xi) - \varepsilon
\vert\xi\vert^2)
\nu (\tau,\xi)\, d\xi\; .\]
For (RVP) the kinetic energy reads
\[ \ekin (\tau) = \int\!\!\!\!\!\int \vert\eta\vert^2 F
(\tau,\xi,\eta)\, d\eta\,d\xi, \]
while for (REP) it reads
\[ \ekin (\tau) =\int \vert\eta\vert^2 (\tau,\xi)\, \nu
(\tau,\xi)\, d\xi . \]
Recalling the remark on the relation between (VP) and (EP) from the
introduction, the second formula can be viewed as a special case
of the first one, and for both
systems we find after a standard computation:
\be \label{eneq}
\frac{d}{d\tau} \Bigl(\ekin (\tau) + \epot (\tau)\Bigr)
= (d-4)\, R^{\frac{d}{2} -1} \dot R \, \ekin (\tau) \quad\mbox{for}\quad d\geq
3\; .\ee
This result is easily achieved as follows. First integrate the rescaled
Vlasov equation  \wrt $\eta$,
\be \label{localcons}
\dtau \nu + \div (\nu \,\eta) = 0,\ee
which is nothing else than the local conservation of mass. Then multiply the rescaled Vlasov equation by
$(|\eta|^2-\varepsilon|\xi|^2)$ and integrate \wrt $\xi$ and $\eta$, and then by
parts,
\[ 
\int\!\!\!\!\!\int (|\eta|^2-\varepsilon|\xi|^2)(\eta\cdot\partial_\xi
F+\varepsilon\xi\cdot\partial_\eta F)\; d\xi d\eta = 0,
\]
so
\[ 
\frac d{d\tau}\!\!\int\!\!\!\!\!\int (|\eta|^2-\varepsilon|\xi|^2)F\; d\xi d\eta 
+2\int\!\!\!\!\!\int \eta\cdot\partial_\xi W F\, d\xi d\eta = (d-4) R^{\frac{d}{2} -1} \dot R \int\!\!\!\!\!\int |\eta|^2 F\, d\xi d\eta\; . 
\]
An integration by parts \wrt $\xi$ gives 
\[ \int\!\!\!\!\!\int \eta\cdot\partial_\xi W \;F\; d\xi d\eta = -\int
W(\partial_\xi\cdot\int \eta\;F\;d\eta)\; d\xi \; , \]
which using \eqn{localcons} gives \eqn{eneq}. Note that \eqn{localcons} is
written for instance in the sense of distributions and that in the above
integration by parts one has to check that no boundary term appears. This 
is true
for $d\geq 3$, but not for $d=2$ as we shall see below.

Recall also that $R^{d/2}\dot R = R^{d/2}dR/dt=d[R(t(\tau))]/d\tau$. Let us
rewrite the energy for the rescaled systems in terms of the original variables:
if we define $P$ and $K$ by 
\[ 
P(t)=\epot (\tau(t))
= R^{d-2} (t) \int \left( U(t,x) - \varepsilon \frac{\vert
x\vert^2}{R^2(t)}
\right)\, \rho (t,x)\, dx, 
\]
and
\[ K (t)=\ekin (\tau(t))
= R^{d-2}(t) \int\!\!\!\!\!\int \biggl| v - \frac{\dot R}{R} x\biggr|^2
f(t,x,v)\, dx\, dv \]
for (VP) or, for (EP)
\[ K (t)=\ekin (\tau(t))
= R^{d-2}(t) \int\
\biggl| u(t,x) - \frac{\dot R}{R} x\biggr|^2
\rho (t,x)\, dx\; , \]
then because of \eqn{eneq}
\be \label{defLyap}
L(t)=K(t)+P(t) 
\ee
is a non-increasing quantity \wrt $t$ for $d=3$,~$4$:
\be \frac {dL}{dt}=(d-4)\frac{\dot R}R K\leq 0\; .\label{Lyap}\ee
Because of the integrations by parts in the intermediate computations, the above
formulas are true only for $d\geq 3$. We will now consider the cases
$d=1$ and $d=2$.

In dimension $d=1$ with $\varepsilon=-1$ (plasma physics case), direct computations
involving the kinetic energy and integral quantities related to the force field
have been used in \cite{BKR} to prove the exponential convergence (in the
rescaled time variable $\tau$) of $F(\tau,\cdot,\cdot)$ towards $F_\infty$ in
$(W^{1,\infty}(\R^2))'$ and of $\partial_\xi W(\tau,\cdot)$ towards $\partial_\xi
W_\infty$ in $L^2(\R)$. The same computation also holds true for the solution of
(EP) if it exists globally in time:

\begin{prop}\label{EP1}
Assume that $d=1$, $\varepsilon=-1$ and consider a global solution 
$(t,x)\mapsto (\rho(t,x),u(t,x))$ of (EP) in $C^1(\R^+\times\R)$ such that 
for any $t>0$, $\rho(t,\cdot)$ has a compact support. Then 
\[ \nu(\tau(t),\xi)=R(t)\rho(t,R(t)\xi)\; ,\quad
\eta(\tau(t),\xi)=\frac 1{\sqrt{R(t)}}\biggl(u(t,R(t)\xi)-\dot{R}(t)\xi\biggr) \]
with $\tau(t)=2\log(1+t)$ and $R(t)=(1+t)^2$ is a solution of (REP) and
converges to $(\nu_\infty^M, 0)$ where $\nu_\infty^M$ is given by Equation
(\ref{asdens}), with $M=\Vert \rho(t,\cdot)\Vert_{L^1(\R)}$: there exists a positive
constant $C$ such that
\[ \Vert\nu(\tau,\cdot)-\nu_\infty^M \Vert_{(W^{1,\infty})'}\leq C\cdot
e^{-\tau}\; , \]
while the electric field 
$
-\partial_\xi W(\tau,\cdot)=\int_{-\infty}^\xi\nu(\tau,\zeta)\; d\zeta 
- \frac 12 \Vert\nu(\tau,\cdot)\Vert_{L^1(\R)}
$
converges in $L^2(\R)$ to 
$-\partial_\xi W_\infty^M$ which is given by Equation (\ref{asfield}):
\[ \Vert\partial_\xi W(\tau,\cdot)-\partial_\xi W_\infty^M \Vert_{L^2(\R)}
\leq 
C\cdot e^{-\tau}\; . \]
In terms of the original variables and with the notation of Section 3, 
this means:
\[ \Vert(1+t)^2\rho(t,(1+t)^2 \cdot)-\nu_\infty^M \Vert_{(W^{1,\infty})'}\leq 
\frac C{(1+t)}\; , \]
\[ \Vert\partial_xU(t,(1+t)^2 \cdot)-\partial_\xi W_\infty^M
\Vert_{L^2(\R)}
\leq\frac C{(1+t)}\; . \]
\end{prop}

Note that in Proposition \ref{EP1}, we made for $R(t)$ the same choice 
as in 
\cite{BKR}, which means that with the notation of Remark \ref{ODE} we 
consider the solution of Equation (\ref{greq}) corresponding to $R_0=1$ and
$\dot{R}_0=2$.
\medskip

The proof essentially follows the same arguments as in \cite{BKR}.
\medskip

In the case $\varepsilon=+1$ (gravitational case), essentially nothing is known
concerning the asymptotic behavior of the solution. If $d=2, 3, 4$ and
$\varepsilon=-1$, the question of identifying the limit of
$F(\tau,\cdot,\cdot)$ or $\nu(\tau,\cdot)$ in the sense of measures as 
$\tau\to\tau_\infty=\int_0^{+\infty} R^{-d/2}(t)\,dt$
(which is finite as soon as $d\geq 3$) is an open question. As already
noted, a natural conjecture would be to identify this limit with
$F_\infty^M$ for the solution of (RVP) and $\nu_\infty^M$ for the solution
of (REP) as in dimension $d=1$. In other terms, the stationary state of
the rescaled equation would be an attractor for the solutions of the
rescaled system in dimension $d>1$. If $d\geq 3$, this is not true in
general.\medskip

\noindent{\bf Counter-examples.} {\sl Consider a solution for which it is
the case and shift the initial data by a constant velocity. Since 
asymptotically the
support of the unscaled solution grows linearly in time,
after rescaling, the shifted solution cannot converge to the stationary
profile. One may then ask the same question in the reference frame of the
center of mass. The following counter-example for (EP) again shows that
for $d\geq 3$, $\varepsilon=-1$, the answer is negative.

Consider in $\R^3$ the solution corresponding to the following initial
data:
\beas 
&&\rho(t=0,x)=3\un_{B^3 (1)}(x)+\un_{B^3 (3)\backslash B^3 (2)}(x)\; ,\\
&&u(t=0,x)=0\quad\mbox{if}\quad \vert x\vert<1\; ,\\
&&u(t=0,x)=x\quad\mbox{if}\quad 2<\vert x\vert<3\; .
\eeas
For any $t>0$, the solution is supported in the union of a centered ball of 
radius $R(t)$ (which obeys to Equation (\ref{greq})) and of a centered
annulus of inner radius $R_1(t)$. A straightforward computation shows that $R$
and $R_1$ satisfy
\[
\ddot{R}=\frac 1{R^2},\; R(0)=1,\; \dot{R}(0)=0,\quad 
\ddot{R}_1=\frac 1{R_1^2},\; R_1(0)=2,\; \dot{R}_1(0)=2\; ,
\]
respectively, and an integration \wrt $t$ gives
\[\dot{R}^2(t)=2-\frac 2{R(t)}<2<4<5-\frac 2{R_1(t)}=\dot{R}_1^2(t)\]
for any $t>0$. As $t\to +\infty$, $\sqrt 2 =\lim_{\to +\infty}\frac{R(t)}t
<\lim_{\to +\infty}\frac{R_1(t)}t=\sqrt 5$ which again forbids the convergence to
the stationary solution after rescaling.\/}
\medskip

In dimension $d=2$ for $\varepsilon=-1$, the situation is different: 
$t\mapsto R(t)$ grows superlinearly (as for $d=1$: see Remark \ref{ODE}, {\sl
(v)\/}), and the question is still open.
\medskip

Consider now the case $d=2$ for (VP). The main difficulty comes from the
integration by parts, and one has to be very careful with the terms involving the
self-consistent potential $U$ since $\nabla U$ essentially decays like $1/|x|$.
Let $(\rho, \rho u) = \int f(t,x,v) (1,v)\; dv$, $(\nu, \nu\eta) = \int
F(\tau, \xi,\eta) (1,\eta)\; d\eta$ and $M=\Vert
f(t,\cdot,\cdot)\Vert_{L^1(\R^2\times\R^2)}$.
\begin{eqnarray} 
\int\frac\xi{|\xi |^2}\cdot (\nu\eta)(\tau,\xi)\; d\xi 
&=& 
\int\frac
{Rx}{|x|^2}\cdot R^2\rho(t,x)(u(t,x)-\frac{\dot R}R x)\; \frac{dx}{R^2}\;
\nonumber\\ 
&=&
-{\dot R}\int\rho(t,x)\; dx+R\int \frac x{|x|^2}
\cdot (\rho u)(t,x)\; dx\nonumber\\ 
&=&
-M\dot R -R\int\frac 1{|x|}\partial_t\rho(t,x)\; dx\nonumber
\end{eqnarray}
using the local conservation of mass $\partial_t\rho+\partial_x(\rho u)=0$.
Similarly,
\begin{eqnarray} 
\int\partial_\xi W\cdot(\int\eta F(\tau,\xi\eta)\; d\eta)\;d\xi 
&=& 
-\frac {M^2}{2\pi}\frac {\dot{R}}R-\int U(t,x)\partial_t\rho(t,x)\;
dx\nonumber\\ 
&=& 
-\frac {M^2}{2\pi}\frac {\dot{R}}R -\frac 12 \frac d{dt}\int
\rho U(t,x)\; dx\nonumber\\ 
&=& 
-\frac {M^2}{2\pi}\frac {\dot{R}}R -\frac 12 \frac d{d\tau}\int
W\eta(\tau,\xi)\; d\xi; \nonumber
\end{eqnarray}
cf.\ \cite{D1} for more details. 
Of course the computations are exactly the
same for (EP). Thus in dimension $d=2$, the definition \eqn{defLyap}
has to be replaced by 
\be L(t)=K(t)+P(t)+\frac {M^2}{2\pi}\log R(t)\; ,\label{defLyap2d}\ee
so that Equation \eqn{Lyap} still holds, see also Remark \ref{direct2d}.
\medskip

Equation (\ref{Lyap}) provides an identity which is a
sharpened form of the Lyapunov functional (also called pseudo-conformal law: see
Section 10 for the relation with the Schr\"odinger equation). A simple form of
this identity had been discovered independently by R.~Illner and G.~Rein, and 
by B.~Perthame, cf.\ \cite{IR,P1}. The improved Lyapunov functional has the 
striking property that it easily provides all the terms that one has to
take into account in the case $d=2$ (see \cite{D1} for (VP)) in a quite 
straightforward manner, while a direct approach was far from being obvious.
\begin{thm}\label{LyaVPp2-4} Assume that $f$ is a solution 
of (VP) with $M=\Vert
f(t,\cdot,\cdot)\Vert$ and that $t\mapsto R(t)$ is the solution of
Equation (\ref{req}) with $R(0)=1$, $\dot{R}(0)=0$.
The function $t\mapsto L(t)$ given by 
\[ 
L(t)=R^{d-2}\int\!\!\!\!\!\int_{\R^d\times\R^d} \biggl| v - \frac{\dot R}{R} 
x\biggr|^2 f\, dv\, dx +
R^{d-2} \int_{\R^d} \left( U(t,x) - \varepsilon \frac{\vert x\vert^2}{R^d}
\right)\, \rho \, dx 
\]
for $d\geq 3$ 
and 
\[ 
L(t)=\int\!\!\!\!\!\int_{\R^2\times\R^2} \biggl| v - \frac{\dot R}{R} 
x\biggr|^2 f\, dv\, dx + \int_{\R^2} \left( U(t,x) - \varepsilon \frac{\vert
x\vert^2}{R^2}
\right)\, \rho \, dx+\frac {M^2}{2\pi }\log R 
\]
if $d=2$ is decreasing for $d=2, 3$, constant for $d=4$, and for any $d\geq 2$
satisfies 
\[ 
\frac {dL}{dt}= (d-4)\, \dot{R}R^{d-3}
\int\!\!\!\!\!\int_{\R^d\times\R^d} \biggl| v - \frac{\dot R}{R} x\biggr|^2
f\, dv\, dx\; . 
\]
Moreover in the plasma physics case $\varepsilon=-1$, $L$ is bounded from below, 
and for $d=2,3$,
\[ 
\int_0^{+\infty} \dot R(s)R^{d-3}(s)
\biggl(\int\!\!\!\!\!\int_{\R^d\times\R^d} 
\biggl| v - \frac{\dot R(s)}{R(s)} x\biggr|^2
f(s,x,v)\, dv\, dx\biggr)\, ds<+\infty\; . 
\]
\end{thm}

\prf $dL/dt$ has already been computed above. For $\varepsilon=-1$, the proof of
the existence of a lower bound is straightforward except maybe for $d=2$. In
that case, the Lyapunov functional is decreasing but might {\sl a priori\/} be
unbounded from below, and we have to estimate it. This can be done using
$M=\ix2{\rho(t,x)}$ and $\ix2{\rho(t,x)\vert x\vert^2}$ by splitting the
integral \wrt $x$ and $y$ into two parts corresponding to $\vert x
-y\vert\leq k$ and $\vert x -y\vert>k$, but a more straightforward approach can
be deduced from Jensen's inequality using the fact that $-\log$ is a convex
function:
\beas 
&&
-\frac 1{2\pi}\int\!\!\!\!\!\int_{\R^2\times\R^2}
\log\vert x-y\vert\, \rho(t,x)\rho(t,y)\; dxdy \\
&& \qquad \qquad \qquad = 
\frac{M^2}{4\pi}\int\!\!\!\!\!\int_{\R^2\times\R^2}
\biggl(-\log\vert x-y\vert^2\biggr)\cdot \rho(t,x)\rho(t,y)\; 
\frac{dxdy}{M^2}\\
&& \qquad \qquad \qquad \geq 
-\frac{M^2}{4\pi}\log\biggl(\int\!\!\!\!\!\int_{\R^2\times\R^2}
\vert x-y\vert^2\cdot \rho(t,x)\rho(t,y)\; \frac{dxdy}{M^2}\biggr)\\
&& \qquad \qquad \qquad \geq 
-\frac{M^2}{4\pi}\log(2I/M)\; ,
\eeas
where $I=\int|x|^2\rho\; dx$, and an optimization on $I>0$ gives
\[ 
\frac{M^2}{2\pi}\log R +\frac I{R^2}-\frac{M^2}{4\pi}\log(2I/M)\geq
\frac {M^2}{4\pi}[1-\log(M/2\pi)]\; ,
\]
which proves the result. \square

\medskip
\begin{rem}\label{Lyapgrav}
For $\varepsilon=+1$ and $d=3$ or $4$, using the Hardy-Littlewood-Sobolev
inequality and classical interpolation identities, one proves that the
self-consistent potential energy term $\Vert\nabla U\Vert^2_{L^2}$ is bounded in
terms of $K$ by
\[ 
\Vert\nabla U\Vert^2_{L^2(\R^d)}\leq C\Vert
f\Vert_{L^1(\R^d\times\R^d)}^{2(1-(d^2-4)/4d)}\Vert
f\Vert_{L^\infty(\R^d\times\R^d)}^{(d-2)/d} K(t)^{(d-2)/2}\; , 
\]
see Section 5 for more details on interpolations.
\end{rem}

\medskip
\begin{rem}\label{direct2d} In dimension $d=2$, for $\varepsilon=-1$, it is probably
easier to compute $dL/dt$ and prove \eqn{Lyap} directly from (VP) using the
identity 
\[ 
\int\!\!\!\!\!\int_{\R^2\times\R^2}(x\cdot v)\; (\partial_xU\cdot \partial_v f)\;
dv\, dx= -\int_{\R^2}(x\cdot \partial_xU)\rho\; dx=\frac {M^2}{4\pi} 
\]
once the equation for $R$ is known, cf.\  \cite{D1}. 
\end{rem}

Note that with the help of \eqn{Requiv} and the results of Theorem
\ref{LyaVPp2-4}, we recover the results of \cite{IR,P1} in dimension
$d=3$ as well as the results of \cite{D1} in dimension $d=2$. Very similar
results of course hold for (EP) since the estimates on the Lyapunov functional in
dimension $d=2,3,4$ are the same.

\begin{thm}\label{LyaEPp2-4} 
Assume that $(\rho,u)$ is a global 
strong solution of (EP) with $M=\Vert\rho(t,\cdot)\Vert_{L^1(\R^d)}$ and that
$t\mapsto R(t)$ is the solution of Equation (\ref{req}) with $R(0)=1$ and
$\dot{R}(0)=0$. The function 
$t\mapsto L(t)$ given by 
\[
L(t)=R^{d-2} \int_{\R^d} \biggl| u - \frac{\dot R}{R} x\biggr|^2\rho\, dx 
+ R^{d-2} \int_{\R^d} \left( U - \varepsilon \frac{\vert x\vert^2}{R^d} \right)\,
\rho \, dx
\]
for $d\geq 3$, and
\[
L(t)= \int_{\R^2} \biggl| u - \frac{\dot R}{R} x\biggr|^2\rho\, dx 
+ \int_{\R^2} \left( U - \varepsilon \frac{\vert x\vert^2}{R^2} \right)\,
\rho \, dx\;+\;\frac {M^2}{2\pi}\log R
\]
for $d=2$, is decreasing for $d=2, 3$, constant for $d=4$, and 
for any $d\geq 2$ satisfies 
\[ 
\frac {dL}{dt}(t)= (d-4)\, \dot{R}R^{d-3}
\int_{\R^d} \biggl| u(t,x) - \frac{\dot R}{R} x\biggr|^2
\rho(t,x)\, dx\; . 
\]
Moreover in the plasma physics case $\varepsilon=-1$, 
$L$ is bounded from below, and for $d=2, 3$
\[ 
\int_0^{+\infty} \dot R(s)R^{d-3}(s)
\biggl(\int_{\R^d} \biggl| u(s,x) - \frac{\dot R (s)}{R(s)} x\biggr|^2
\rho(s,x)\, dx\biggr)\, ds <+\infty\; . 
\]
\end{thm}
\bigskip

The case $d=4$ appears to be the limit case to which the above 
method for finding Lyapunov functionals in the plasma physics case applies since
for $d\geq 5$, $t\mapsto L(t)$ is increasing. However for $d\geq 4$, we may write
\[ \frac {dL}{dt}\leq (d-4)\frac{\dot{R}}R L\, , \]
\[ \frac {d}{dt}\biggl(R^2\frac L{R^{d-2}}\biggr)\leq 0\, , \]
and thus obtain 
\[ 
\int\!\!\!\!\!\int_{\R^d\times\R^d} \biggl| v - \frac{\dot R}{R} x\biggr|^2
f(t,x,v)\, dv\, dx=O\left(R^{-2}\right)=O\left(t^{-2}\right) 
\]
since for $d\geq 3$ all the quantities involved in $L(t)$ are nonnegative
and $R(t)\sim t$ as $t\to +\infty$.

In this last part of Section 4, we will derive the Lyapunov
functionals in another way, not because of the case $d>4$ (which is of minor 
interest for (EP) or (VP) in itself), but because the method is simpler 
and will be applied to other systems in Sections 6--8. We assume that
$\varepsilon=-1$ in the rest of this section.

We may indeed notice that all the quantities we have been taking into account
are integrated in the $x$ variable, so that the change of variable
$\xi(t,x)=x/R(t)$ does not play any role in the estimates.
Let us first consider the Vlasov-Poisson system (VP). According to the above
remark, we may use the change of variables
\[ 
\eta(t,x,v)=v-\frac{\dot{R}}Rx\; ,\quad f(t,x,v)=F(t,x,\eta) 
\]
so that $F$ solves the rescaled system (R$'$VP):
\[
\partial_tF+\eta\cdot\partial_xF - \frac{\ddot{R}}{R}
x\cdot\partial_\eta F -\partial_xU(t,x)\cdot\partial_\eta F +
\frac{\dot{R}}{R}\bigl[\partial_x(xF)-\partial_\eta(\eta F)\bigr]=0\; ,
\]
\[
-\partial_xU(t,x)=\frac x{\vert S^{d-1}\vert\, \vert x\vert^d}
*\int_{R^d}F(t,x,\eta)\; d\eta\; .
\]
As for (RVP), we may compute the energy:
\[ 
E(t)=\ixed{\biggl(\vert\eta\vert^2+\frac{\ddot{R}}{R} \vert
x\vert^2+U\biggr)\,F}
\]
if $d\geq 3$ and 
\[ 
E(t)=\ixed{\biggl(\vert\eta\vert^2+\frac{\ddot{R}}{R} \vert
x\vert^2+U\biggr)\,F}
+ \frac{M^2}{2\pi}\frac{\dot{R}}R
\]
for $d=2$. This energy is a decaying function of $t$ : 
for any $d\geq 2$,
\beas 
\frac{dE}{dt}(t)
&=&
- (d-2)\frac{\dot{R}}{R}\ixed{FU} \\ 
&& 
{}+ \ixed{\biggl[ \biggl( \frac d{dt}(\frac{\ddot{R}}{R})
+2\frac{\ddot{R}}{R}\frac{\dot{R}}{R}\biggr)
\vert x\vert^2
- 2\frac{\dot{R}}{R}\vert\eta\vert^2\biggr]\,F }\; .
\eeas
We may now define $L(t)=B(t)E(t)$. 
For any $d\geq 3$,
\begin{eqnarray}
\frac{dL}{dt}
&=&
\biggl(\dot{B}-(d-2)\frac{\dot{R}}{R}B\biggr)
\ixd{\vert \nabla U\vert^2}\nonumber\\
&&
{}+ \biggl(\dot{B} - 2\frac{\dot{R}}{R}B\biggr)
\ixed{F(t,x,\eta)\vert\eta\vert^2}\\
&&
{}+ \biggl(\dot{B}\frac{\ddot{R}}{R} + 
\biggl(\frac d{dt}(\frac{\ddot{R}}{R}) + 2\frac{\ddot{R}}{R}\frac{\dot{R}}{R}\biggr)\,B\biggr)
\ixed{F(t,x,\eta)\vert x\vert^2}\; ,\nonumber\label{Ldirect}
\end{eqnarray}
while for $d=2$, 
\begin{eqnarray}
\frac{dL}{dt}
&=&
\biggl(\dot{B}\log R-B\frac{\dot{R}}R\biggr) M^2\nonumber\\
&&
{}+ \biggl(\dot{B} - 2\frac{\dot{R}}{R}B\biggr)
\ixe2{F(t,x,\eta)\vert\eta\vert^2}\\
&&
{}+ \biggl(\dot{B}\frac{\ddot{R}}{R} +
\biggl(\frac d{dt}(\frac{\ddot{R}}{R}) +
2\frac{\ddot{R}}{R}\frac{\dot{R}}{R}\biggr)\,B\biggr)
\ixe2{F(t,x,\eta)\vert x\vert^2}\; .\nonumber\label{Ldirect2}
\end{eqnarray}
For $d\geq 3$, the following conditions are sufficient for $L$ to be
nonincreasing:
\begin{description}
\item{1)} 
$B(t)=R(t)^{d-2}$, which implies
$\dot{B}-(d-2)\, B\, \dot{R}/R\leq 0$,
\item{2)} 
$d\leq 4$, which implies
$\dot{B} - B\,\dot{R}/R = -(4-d)\, B\, \dot{R}/R \leq 0$ ,
\item{3)} 
$\ddot{R}=R^p\; ,\quad R(0)=1\; ,\quad\dot{R}(0)=0$ with $p\leq -(d-1)$, which 
implies
$(\frac d{dt}(\frac{\ddot{R}}{R})+2\frac{\ddot{R}}{R}\frac{\dot{R}}{R})B +
\dot{B}\frac{\ddot{R}}{R}\leq 0$,
\end{description}
and we recover the results of Theorem (\ref{LyaVPp2-4}) for $d=3, 4$;
for $d=2$ take $B=1$. 

\begin{rem}\label{dgeq5} If $d\geq 2$ (including the case $d\geq 4$), we may
choose
$B=R^{d-2-\theta}$, $\theta\geq \max(0,d-4)$, and $R$ solving the
equation $\ddot{R}=R^p$,
$R(0)=1$, $\dot{R}(0)=0$ for some $p\leq\theta-(d-1)$
without any further restriction on $d$. Note that for
$d\geq 4$, $p<-1$ and $\theta <d-2$, one recovers the
estimate one would have for the free transport $\dt f + v\cdot \dx f =0$, 
since in that case $f(t,x,v)=f_0(x-vt,v)$ and
\[ \ixv{f(t,x,v)\vert x-vt\vert^2} = \ixv{f_0(x,v)\vert x\vert^2}\, . \]
For the consequences on the dispersion rate, see Section 5.
\end{rem}
\medskip

An analogous method also works for the Euler-Poisson system (EP). If
we shift the velocity $u(t,x)$ by an unknown "bulk" velocity 
$\frac{\dot{R}}R x$, so that $\eta(t,x)=u(t,x)-\frac{\dot{R}}Rx$,
then $(\rho(t,x),\eta(t,x))$ solves the system (R$'$EP):
\[
\partial_t\rho+\partial_x(\rho (\eta+\frac{\dot{R}}Rx))=0\; ,
\]
\[
\partial_t\eta + \frac d{dt} (\frac{\dot{R}}R)x
+ ((\eta+\frac{\dot{R}}Rx)\cdot\partial_x)\eta
+\frac{\dot{R}}R(\eta+\frac{\dot{R}}Rx)=-\partial_xU(t,x)\; ,
\]
\[
-\partial_xU(t,x)=\frac x{\vert S^{d-1}\vert\, \vert x\vert^d}*\rho\; .
\]
As for (RVP), we may consider the energy:
\[ 
E(t)=\ixd{\biggl(\vert\eta(t,x)\vert^2 + \frac{\ddot{R}}{R}
\vert x\vert^2 + U(t,x)\biggr)\, \rho(t,x)} 
\]
for $d\geq 3$ (if $d=2$, one has to add the term
$\frac{M^2}{2\pi}\frac {\dot{R}}R$). 
Again $t\mapsto E(t)$ is decaying: for any $d\geq 2$,
\[ 
\frac{dE}{dt}=-(d-2)\frac{\dot{R}}{R}\ixd{\rho U} +
\ixd{\biggl[ \biggl(\frac
d{dt}(\frac{\ddot{R}}{R})+2\frac{\ddot{R}}{R}\frac{\dot{R}}{R}\biggr)|x|^2
- 2\frac{\dot{R}}{R}\vert\eta\vert^2\biggr]\, \rho }\, . \]
As for (VP), we may also define $L(t)=B(t)E(t)$, and the rest of the 
discussion is exactly the same.\medskip

The method for finding a Lyapunov functional can be summarized as follows:
first change the velocity variable by subtracting a velocity $\frac{\dot{R}}R
x$ for some increasing function $R$, then compute the energy associated to the
new equation and finally choose the Lyapunov functional to be 
$L(t)=B(t) E(t)$ 
where $B(t)$ is the function of $t$ which has the maximal growth in 
order that $L(t)$ is still a decaying function of $t$ and corresponds to 
a function $t\mapsto R(t)$ solving an adequate ordinary differential equation
which takes the nonlinearity into account and has to be chosen well.
Of course, one way to find an equation for $R$ is to apply the method
of the time-dependent rescalings of the beginning of this section.
This method is sufficient to extract the asymptotic rate of decay of
the relevant quantities, as we shall see later in several other cases,
cf.\ Sections 6--8.

\section{Asymptotic behaviour, dispersion}
\setcounter{equation}{0}

An estimate of the rate of dispersion of a solution $f$ of the Vlasov-Poisson 
system (VP) in the plasma physics case $\varepsilon=-1$
is given by the interpolation of $\rho(t,x)=\int_{\R^d} f(t,x,v)\; dv$ 
between the $L^\infty$-norm of $f$, which is preserved for strong solutions,
and the momentum 
$\int\!\!\!\!\!\int_{\R^d\times\R^d} {f\vert v-\frac
xt\vert^2}\; dxdv\sim \int\!\!\!\!\!\int_{\R^d\times\R^d}{f\vert v-\frac 
{\dot{R}}R x\vert^2}\; dxdv$ as $t\to +\infty$: there exists a constant
$C=C(d)>0$ such that
\be \label{interpolation}
\biggl|\biggl|\int_{\R^d} f\; dv\biggr|\biggr|_{L^{\frac d{d+2}}(\R^d)}
\leq
C\cdot\Vert f\Vert_{L^\infty(\R^d\times\R^d)}^{\frac
2{d+2}}\cdot\biggl(\int\!\!\!\!\!\int_{\R^d\times\R^d}{f\biggl|
v-\frac{\dot{R}}Rx\biggr|^2}\; dxdv\biggr)^{\frac d{d+2}};
\ee
for a systematic study of these interpolation inequalities
see \cite{D1} and references therein.

The asymptotic form of the Lyapunov functional was given
in \cite{IR,P1} for the case $d=3$ and in
\cite{D1} for the case $d=2$. Using $R(t)$, we remove the
difficulty due to the singularity at $t=0$ and recover the
known results. The use of the decay term of the
Lyapunov functional allows us to prove that the decay is not
optimal.

\begin{prop}\label{VPdisp} Assume that $f$ is a strong solution of (VP)
in the plasma physics case $\varepsilon=-1$ and $t\mapsto R(t)$ is the solution 
of $\ddot R = R^{1-d}$ with $R(0)=1$ and $\dot{R}(0)=0$. Then $f$
obeys to the following Strichartz type estimate: if $d=2$ or $3$,
\be\label{int3}
\int_0^{+\infty}R^{d-3}(t)\dot{R}(t)\biggl(\int\!\!\!\!\!\int_{\R^d\times\R^d}
{f(t,x,v)\biggl| v-\frac{\dot{R}(t)}{R(t)}x\biggr|^2}\; dxdv\biggr)\; dt
\leq C\; ,
\ee
and for $d=3, 4$, we have the following dispersion estimate
\be\label{decay34}
\Vert \rho(t,\cdot)\Vert_{L^{\frac{d+2}d}(\R^d)}
\leq
C \, R(t)^{-d \frac {d-2}{d+2}}\sim t^{-d \frac {d-2}{d+2}}\; .
\ee
Here $C$ denotes various positive constants which depend only on $d$ 
and $f_0$, and $L$ is the Lyapunov functional of Theorem \ref{LyaVPp2-4}.
If $d=2$, 
\begin{equation}\label{decay2}
\liminf_{t\to+\infty}\; \Vert \rho(t)\Vert_{L^2(\R^2)} =0\; .
\end{equation}
\end{prop}

The proof follows from Theorem \ref{LyaVPp2-4} and the interpolation
identity (\ref{interpolation}) given above. The decay of $\rho(t,\cdot)$ in
$L^2(\R^2)$ is given by the decay term of the Lyapunov functional; see
Remark \ref{opt} below. Estimate \eqn{int3} for $d=2$ has
been improved compared to \cite{D1}.

\begin{rem} \label{opt} The decay given in Proposition \ref{VPdisp} is 
not optimal. Consider indeed a function 
$t\mapsto h(t)$ such that $h\geq 1$, $\lim_{t\to
+\infty}h(t)=+\infty$ and 
\[ 
\int_0^{+\infty}\frac
{ds}{s\, h(s)}=+\infty\quad\mbox{if\ }d=2\; ,\quad 
\int_0^{+\infty}\frac {ds}{h(s)}=+\infty
\quad\mbox{if\ }d=3\; .
\]
For instance, one may take for $t>0$
\[ 
h(t)=\log(t+2)\quad\mbox{if}\ d=2\quad\mbox{and}\quad 
h(t)=t\, \log(t+2) \quad\mbox{if}\ d=3\; .
\]
The bounds (\ref{int3}) immediately provide for $d=2, 3$
\[ 
\liminf_{t\to +\infty} \, h(t)\,
\int\!\!\!\!\!\int_{\R^d\times\R^d}{f(t,x,v)
\biggl| v - \frac{\dot{R}(t)}{R(t)}x\biggr|^2}\; dxdv = 0, 
\]
and as a consequence, the decays in (\ref{decay34}) and (\ref{decay2}) are
not optimal.
\end{rem}

Similar results can of course be obtained for any $d\geq 5$, using Remark
\ref{dgeq5}. We may notice that the decay in (\ref{decay34}) for $d\geq 4$
is the one which is obtained for the free transport equation when 
considering the second moment in $x-vt$.

\begin{rem}\label{equiv}
The Lyapunov functionals given
in \cite{D1,IR,P1} correspond to the
asymptotic form of $R(t)$ as $t\to +\infty$. The fact that this
asymptotic form also gives a Lyapunov functional is easily explained
by the scaling invariance of the equation (see Remark \ref{ODE}): 
if one replaces $R(t)$ by $t$ for $d=3, 4$ or $t\sqrt{\log t}$ for 
$d=2$ in the expression of the Lyapunov functional $L(t)$ of Theorem 
\ref{LyaVPp2-4}, $L(t)$ would still be a Lyapunov functional (consider its limit
as $\lambda\to +\infty$).
\end{rem}
\medskip

Similar results for the pressureless Euler-Poisson system (EP) also hold 
except that no direct interpolation can be used. The decay only holds in 
a weak norm defined as follows (assume here that $d\geq 3$): let us consider
the space
\[ 
{\cal D}^{1,2}(\R^d)=\left\{\phi\in 
L^{\frac{2d}{d-2}}(\R^d)\; :\; \nabla\phi\in L^2(\R^d)\right\} . 
\]
and define on its dual space the norm
\[ 
\vert\vert\vert \rho \vert\vert\vert=
\Vert\rho\Vert_{({\cal D}^{1,2}(\R^d))'}=\sup\biggl\{\int\rho\phi\, dx
\mid \phi\in {\cal D}^{1,2}(\R^d),\ 
\Vert\nabla\phi\Vert_{L^2(\R^d)}\leq 1\biggr\}\; . 
\]
If $U\in{\cal D}^{1,2}(\R^d)$ is such that
$-\Delta U=\rho$, then
$\vert\vert\vert \rho \vert\vert\vert
\leq \Vert\nabla U\Vert_{L^2(\R^d)}$.
Using the same notation as in Section 4, if $d=3, 4$, there exists a
positive constant $C$ such that
\[ 
\ixd{\rho(t,x)\vert\eta(t,x)\vert^2} =
R(t)^{d-2}\ixd{\rho(t,x)
\biggl| u(t,x)-\frac{\dot{R}}Rx \biggr|^2}\leq C 
\]
and
\[ 
R(t)^{d-2}\ixd{\rho(t,x)U(t,x)} =
R(t)^{d-2}\ixd{\vert \nabla U(t,x)\vert^2}\leq C\; . 
\]
The last inequality can be reinterpreted as an estimate on the 
weak norm $\vert\vert\vert\cdot \vert\vert\vert$ of $\rho(t,\cdot)$.

\begin{prop} 
Assume that $(\rho,u)$ is a $C^1$ solution on $\R^+\times\R^d$
of (EP) in the plasma physics case $\varepsilon=-1$ and $t\mapsto R(t)$ is 
the solution of $\ddot R + \varepsilon R^{1-d}=0$ with $R(0)=1$ and $\dot{R}(0)=0$. 
Then $(\rho,u)$ obeys to the following Strichartz type estimate for $d=2$, $3$:
\[ 
\int_0^{+\infty}R^{d-3}(t)\dot{R}(t)
\biggl(\ixd{\rho(t,x)\vert
u(t,x) -\frac{\dot{R}(t)}{R(t)}x\vert^2}\biggr)\; dt\; 
<\; +\infty\; . 
\]
Moreover, if $d=3, 4$, then
\[ 
\limsup_{t\to +\infty}\,
t^{\frac d2-1}\vert\vert\vert
\rho(t,\cdot)\vert\vert\vert=\limsup_{t\to +\infty}\,
R(t)^{\frac d2-1}\Vert
\nabla U(t,\cdot)\Vert_{L^2(\R^d)}< +\infty\; . 
\] 
\end{prop}

\begin{rem}
If $d=2$, we cannot use the $({\cal D}^{1,2}(\R^d))'$-norm as in the case 
$d\geq 3$, but the following estimates for the solutions in
the plasma physics case $\varepsilon=-1$ of the 
pressureless Euler-Poisson system (EP) hold:
\[
\lim_{t\to +\infty}\frac 1{\log R(t)}\ix2{\rho(t,x)U(t,x)}=-\frac
{M^2}{2\pi}\; ,
\]
\[
\lim_{t\to +\infty}\frac 1{\log R(t)}\ix2{\rho(t,x)\vert
u(t,x)\vert^2}=\frac {M^2}{2\pi}\; ,
\]
\[
\lim_{t\to +\infty}\frac 1{t^2}\ix2{\rho(t,x)\vert
x\vert^2}=\frac {M^2}{2\pi}\; .
\]
These estimates are easily deduced from the conservation of
the energy, the expression of the Lyapunov functional $L(t)$ and the
estimate given in the proof of Theorem \ref{LyaVPp2-4}.
\end{rem}

Maybe more interesting is the observation (see \cite{D1}) that for
$d=2$, which is the limit case for dispersion results, the dispersion
estimate gives a lower bound for the growth of the support of a
solution corresponding to a compactly supported initial datum:

\begin{cor} Consider for $d=2$ solutions of (VP) or
(EP) corresponding to compactly supported initial data.
Assume that $r(t)$ is the minimal radius of the balls containing
the support of $\rho(t,\cdot)$.  
Then there exists a constant $C>0$ such that 
\[ r(t)\geq C\, R(t)\quad\mbox{as}\quad t\to +\infty\; . \]
\end{cor}

\prf As in \cite{D1} one may simply notice that 
\[ 
\frac {M^2}{2\pi}\biggl(\log R(t) -\log(2r(t))\biggr)\leq L(t)\leq L(0)
\; .\qquad\Box
\]

\section{The $2$-dimensional symmetric Vlasov-Pois\-son system with 
an external magnetic field}
\setcounter{equation}{0}

In dimension $d=2$, we may consider the following system (VPM)
\beas
&&\partial_tf+v\cdot\partial_xf+
\biggl(-\partial_xU(t,x)+B_0v^\perp\biggr)\cdot\partial_vf=0\\
&&-\partial_xU(t,x)=\frac x{2\pi\,
\vert x\vert^2}*\int_{\R^2}f(t,x,v)\; dv
\eeas
corresponding to a system of particles with a
self-interaction through electrostatic forces, in the
presence of an external constant magnetic field $B_0$. Here
we use the notation
\[ 
\left(\begin{array}{c}v_1\\v_2\end{array}\right)^\perp 
=
\left(\begin{array}{c}-v_2\\v_1\end{array}\right)\, . 
\]
For the linear system without self-consistent electrostatic forces, 
all the characteristics are circles and a solution 
with an initially compact support will remain supported in
a fixed compact set for all time. With a self-consistent
Poisson term, the situation is radically different since we get 
the same estimates as for the Vlasov-Poisson system without a 
magnetic field.

We may indeed shift the velocity variable
$\eta(t,x,v)=v-\frac{\dot{R}}R x$,
and the new distribution function $f(t,x,v)=F(t,x,\eta)$
obeys to the system
\beas
\partial_t F+\eta\cdot\partial_x F 
- 
\frac{\ddot{R}}{R}
x\cdot\partial_\eta F 
+\biggl(-\partial_x U(t,x) + B_0 (\eta^\perp
+\frac{\dot{R}}{R}x^\perp)\biggr)\;
\cdot\partial_\eta F \ &&\\
+ \frac{\dot{R}}{R}\biggl(\partial_x(xF)-\partial_\eta(\eta F)\biggr)
&=&
0\; ,
\eeas
\[
-\partial_xU(t,x)=\frac x{2\pi\, \vert
x\vert^2}*\int_{\R^2} F(t,x,\eta)\; d\eta\; .
\]
Exactly as in Section 4, the energy is
\[ 
E(t)=\frac {M^2}{2\pi}\log R+\int\!\!\!\!\!\int_{\R^2\times\R^2}
F(t,x,\eta)\biggl(\vert\eta\vert^2 +
\frac{\ddot{R}}{R}\vert x\vert^2 + U(t,x)\biggr)\,d\eta\, dx 
\]
and decays according to
\[ 
\frac{dE}{dt}= - \ixed{F(t,x,\eta)\biggl[2\frac{\dot{R}}{R}\vert\eta\vert^2 +
\biggl(\frac d{dt}(\frac{\ddot{R}}{R})+2\frac{\ddot{R}}{R}\frac{\dot{R}}{R}\biggr)
\vert x\vert^2 + 2\frac{\dot{R}}R(x\cdot\eta^\perp)\biggr]}\; . 
\]
If $f$ is radially symmetric, {\sl i.e.\/} depends only on $t$,
$\vert x\vert$, $(x\cdot v)$ and $\vert x\vert^2\vert
v\vert^2-(x\cdot v)^2$, then the analogous property holds for
$F$: $F$ only depends on $t$,
$\vert x\vert$, $(x\cdot \eta)$ and $\vert x\vert^2\vert
\eta\vert^2-(x\cdot \eta)^2$, and
\[ 
\int_{\R^2}\biggl(x^\perp\cdot\int_{R^2} \eta F(t,x,\eta)\;
d\eta\biggr)\,dx =
\int\!\!\!\!\!\int_{\R^2\times\R^2}
(x\cdot\eta^\perp)F(t,x,\eta)=0\; . 
\]
The system has the same Lyapunov functional as (VP), 
and we obtain the same dispersion
results as for the Vlasov-Poisson system:

\begin{prop} Let $d=2$. Assume that $f$ is a solution 
of (VPM) and that $t\mapsto R(t)$ is the solution of $\ddot R = \frac 1R$
with $R(0)=1$, $\dot{R}(0)=0$. The function $t\mapsto L(t)$ given by 
\[ 
L(t)=\frac {M^2}{2\pi}\log R + \int\!\!\!\!\!\int_{\R^2\times\R^2} 
\biggl| v - \frac{\dot R}{R} x\biggr|^2 
f(t,x,v)\, dv\, dx +\int_{\R^2} \left(U +
\frac{\vert x\vert^2}{R^2}
\right)\, \rho\, dx 
\]
is decreasing, bounded from below and satisfies
\[ 
\frac {dL}{dt}= -2\, \frac{\dot R}R
\int\!\!\!\!\!\int_{\R^2\times\R^2} \biggl| v - \frac{\dot R}{R} x\biggr|^2
f(t,x,v)\; dv\, dx\; . 
\]
Moreover 
\[ 
\int_0^{+\infty} \frac{\dot R(s)}{R(s)}
\biggl(\int\!\!\!\!\!\int_{\R^2\times\R^2}
\biggl| v - \frac{\dot R(s)}{R(s)} x\biggr|^2
f(s,x,v)\; dv\, dx\biggr)\; ds\; <\; +\infty 
\]
and
\[ 
\liminf_{t\to+\infty}\; \Vert \rho(t) \Vert_{L^2(\R^2)} =0\; . 
\] 
\end{prop}

\section{The isentropic Euler system for perfect gases}
\setcounter{equation}{0}

As another example, which does not belong to the field of kinetic equations, 
we consider the isentropic Euler system (IE)
for perfect gases (for $\gamma>1$)
\[
\partial_t\rho+\partial_x(\rho u)=0\; ,
\]
\[
\partial_t u+(u\cdot\partial_x)u=-\partial_x p\; ,
\]
\[
p=\rho^{\gamma-1}\; .
\]
The method goes exactly as for the pressureless Euler-Poisson system (here we
use the second method of Section 4): the rescaled system (R$'$IE)
given by
$\eta(t,x)=u(t,x)-\frac{\dot{R}}Rx$ is 
\[
\partial_t\rho+\partial_x(\rho (\eta+\frac{\dot{R}}Rx))=0\; ,
\]
\[
\partial_t\eta + \frac d{dt} (\frac{\dot{R}}R)x
+ ((\eta+\frac{\dot{R}}Rx)\cdot\partial_x)\eta
+\frac{\dot{R}}R(\eta+\frac{\dot{R}}Rx)=-\partial_x\rho^{\gamma-1}\; .
\]
The last equation can be rewritten as
\[ \partial_t\eta + \eta\cdot\partial_x\eta + \frac{\dot{R}}R
x\cdot\partial_x\eta + \frac{\ddot{R}}Rx + \frac{\dot{R}}R\eta =
-\partial_x\rho^{\gamma-1}\; .\]
If we define the energy by
\[ E(t)=\ixd{\rho(t,x)\vert\eta(t,x)\vert^2}+\frac{\ddot{R}}{R}\ixd{\rho(t,x)\vert
x\vert^2}+\frac 2\gamma\ixd{\rho^\gamma(t,x)}\; , \]
a Lyapunov functional is easily exhibited by considering
$L(t)=B(t)E(t)$. The energy is indeed decreasing:
\[ 
\frac{dE}{dt}=\ixd{\rho\biggl[ -
2\frac{\dot{R}}{R}\vert\eta\vert^2 + 
\biggl(\frac
d{dt}(\frac{\ddot{R}}{R})+2\frac{\ddot{R}}{R}\frac{\dot{R}}{R}\biggr)
\vert x\vert^2 
- 2d\frac {\gamma-1}\gamma\frac{\dot{R}}R \rho^{\gamma-1}\biggr]} 
\]
so that
\beas
\frac{dL}{dt}
&=&
\biggl(\dot{B} - 2\frac{\dot{R}}{R}B\biggr)
\ixd{\rho(t,x)\vert\eta(t,x)\vert^2}\\
&& + \biggl[\dot{B}\frac{\ddot{R}}{R} +\biggl(\frac
d{dt}(\frac{\ddot{R}}{R})+2\frac{\ddot{R}}{R}\frac{\dot{R}}{R}\biggr)B\biggr]
\ixd{\rho(t,x)\vert x\vert^2}\\
&&+ \frac
2\gamma \biggl(\dot{B}-(\gamma-1)\,d\,B\,\frac{\dot{R}}R\biggr)
\ixd{\rho^\gamma(t,x)} ,
\eeas
and sufficient conditions for $L$ to be decreasing are therefore given by:
\begin{description}
\item{1)} $B=R^q$ with $q\leq \min(2,(\gamma-1)d)$, which implies
$\dot{B} - 2\frac{\dot{R}}{R}B\leq 0$ and $\dot{B}-
d\,(\gamma-1)\, B\, \dot{R}/R\leq 0$.
\item{2)} $\ddot{R}=R^p$ with $p\leq -(q+1)$, which implies
$\dot{B} \ddot{R}/R +(\frac
d{dt}(\frac{\ddot{R}}{R})+2\frac{\ddot{R}}{R}\frac{\dot{R}}{R})B\leq 0$.
\end{description}
It turns out that these dispersion relations (or at least their asymptotic form)
are already known and have been used for 
the Navier-Stokes equation by J.-Y. Chemin in 
\cite{Ch}, and by D. Serre in \cite{GS,Se}
and B. Perthame in \cite{P2}. One of the interests of 
these estimates is that one may use them as an {\sl a priori\/} estimate to
control the behaviour for large times and build a global (in time)
solution to the Cauchy problem. An equivalent remark (see
\cite{Se}) is that it is possible to build a solution by a
fixed-point method for a finite time (this is not in contradiction
with T. Sideris' results \cite {Si} on non-existence, if the initial
data is small in the correct sense) and that one may choose the
rescaling $t\mapsto R(t)$ such that (for the complete rescaling as defined in
Section 3 of course) the evolution \wrt the rescaled time holds only on a
finite time interval
$0\leq\tau<\tau_\infty=\int_0^{+\infty}A^{-2}(t)\,dt$.
However, we are here interested only in the dispersion relations which were 
easily obtained by the mean of the second method of Section 4. These 
dispersion relations can be summarized as follows:

\begin{prop} 
If $(\rho, u)$ is a global classical solution of (IE) 
with $\gamma >1$, then it satisfies the following dispersion relation
\[ 
\frac{d}{dt}\biggl(R^q\ixd{\rho
\biggl| u-\frac{\dot R} R x\biggr|^2} +
\frac{1}{R^2}\ixd{\rho\vert x\vert^2} + \frac 2\gamma R^q\ixd{\rho^\gamma}\biggr) \leq 0 
\] 
with $q=\min(2,(\gamma-1)d)$ and $t\mapsto
R(t)$ such that $\ddot{R}=R^{-(q+1)}$, $R(0)=R_0>0$, $\dot{R}(0)=0$.
\end{prop}

\section{Wigner and Schr\"odinger equations}
\setcounter{equation}{0}

The relation between the Schr\"odinger equation, the Wigner equation 
and the Vlasov equation is now quite well understood. It has been the
subject of a considerable number of papers in the recent years: we 
mention \cite{GMMP,I} as some of the most recent ones, and also 
\cite{LPa,MM} for the limit of the Schr\"odinger-Poisson to the 
Vlasov-Poisson system. Historically, the dispersion relations have been 
studied for the Schr\"odinger equation first and then adapted to the 
corresponding kinetic equation \cite{DF1,DF2,IZL}. The analysis of the 
dispersion relations in the kinetic framework came only after, but 
now seems to provide powerful tools to build new dispersion identities, 
cf.\ \cite{GMP}.

Consider the Schr\"odinger equation
\[ 
i\hbar\partial_t\psi=-\frac 12 {\hbar}^2\Delta \psi+ V\psi\; . 
\]
If $w$ is the Wigner transform of $\psi$,
\[ 
w(t,x,v)=\int_{\R^d}e^{-ivy}\; \overline{\psi}(t,x+\frac\hbar{2}y)
\psi(t,x-\frac\hbar{2}y)\; dy\; , 
\]
it has to satisfy the Wigner equation
\[
\partial_tw+v\cdot\partial_xw-\frac i\hbar\Theta(V)w=0
\]
where the pseudo-differential operator $\Theta(V)$ is defined by
\[ 
\Theta(V)f(x,v)=\frac1{(2\pi)^d}\int_{\R^d}e^{-ivy}\biggl[V({\scriptstyle
x+\frac\hbar{2}y}) -V({\scriptstyle x-\frac\hbar{2}y})\biggr]\cdot
\biggl(\int_{\R^d}e^{+iy\xi}f({\scriptstyle x,\xi})\; d\xi\biggr)\; dy\; . 
\]
In the semi-classical limit $\hbar\to 0+$, the operator
$\Theta(V)$ is formally expected to converge to $-\partial_xV\cdot\partial_v$,
and it is the purpose of many papers  to
justify this limit, cf.\ \cite{GMMP,I,LPa,MM}. 

In this section we will only derive some dispersion identities according 
to the technique developed at the end of Section 4 and give some easy 
consequences of these estimates.
\medskip

We shall consider three cases:
\begin{description}
\item{The linear case {\bf (L)}}: $V$ is a
given fixed nonnegative potential which does not depend on $t$ and decays as 
$\vert x\vert\to +\infty$. We will not go further into this case since the
dispersion properties would depend on the local properties of $V$ and
$x\cdot\partial_x V$, but the computations are essentially the same as for the
other cases up to Equation \eqn{LyapW}.
\item{The Poisson case {\bf (P)}}: In case of
the Schr\"odinger equation $V$ obeys
\[ 
-\Delta V=\vert\psi\vert^2, 
\] 
and in case of the Wigner equation $V$ obeys 
\[ 
-\Delta V=\int_{\R^d}w(t,x,v)\; dv; 
\] 
we consider only the electrostatic case. We shall
state a result on the Wigner and the Schr\"odinger formulations of the problem,
which clearly proves that this case can be handled in full generality with our
methods. The estimates are slightly improved in dimension $d=3$ and can
obviously be generalized to any dimension $d\geq 4$. The results are new for
$d=2$.
\item{The nonlinear case {\bf (NL)}}: $V$ is given by a
power law
\[ 
V=\vert\psi\vert^{p-1}\; .
\] 
In that case we consider only the nonlinear Schr\"odinger equation (NLS)
\[ 
i\hbar\partial_t\psi=-\frac 12 {\hbar}^2\Delta \psi - \varepsilon
\vert\psi\vert^{p-1}\psi 
\]
(in the following, we shall only study the defocusing case $\varepsilon=-1$).
This case is mentioned here to make the link with the pseudo-conformal methods
and to recover the pseudo-conformal law, which has been  studied
extensively. 
\end{description}

\subsection{Wigner equation}

For the Wigner equation, we introduce as for the Vlasov-Poisson system 
the new velocity variable $\eta(t,x,v)=v-\frac{\dot{R}}Rx$ and exactly as for
the Vlasov-Poisson system, $F(t,x,\eta)=w(t,x,v)$ solves the rescaled Wigner
equation (R$'$W):
\[ 
\partial_tF+\eta\cdot\partial_xF - \frac{\ddot{R}}{R}
x\cdot\partial_\eta F - \frac i\hbar\Theta(V)F +
\frac{\dot{R}}{R}
\biggl(\partial_x(xF)-\partial_\eta(\eta F)\biggr)=0\; . 
\]
Again as for (R$'$VP), we compute the energy
\[ 
E(t)=\ixed{F\biggl(\vert\eta\vert^2+\frac{\ddot{R}}{R}\vert
x\vert^2+\alpha V\biggr)}
\]
if $d\geq 3$, and
\[ 
E(t)=\frac{M^2}{2\pi}\log R+ \ixed{F\biggl(\vert\eta\vert^2+
\frac{\ddot{R}}{R}\vert x\vert^2+\alpha V\biggr)}
\]
if $d=2$ in case (P). Here $\alpha$ is a coefficient which 
takes different values
according to the case we consider: $\alpha=2$, $1$, and $\frac 2{p+1}$
in case (L), (P) and (NL) respectively. The same computation as before provides
\beas
\frac{dE}{dt} 
&=&\ixed{\!\!\!F(t,x,\eta)\biggl[-
2\frac{\dot{R}}{R}\vert\eta\vert^2 + 
\biggl(\frac d{dt}(\frac{\ddot{R}}{R}) + 2\frac{\ddot{R}}{R}\frac{\dot{R}}{R}\biggr) 
\vert x\vert^2\biggr]}\\
&&
{}-\alpha\frac{\dot{R}}{R}\ixed{\!\!\! V(t,x)\partial_x(xF(t,x,\eta))}\eeas
(for $d\geq 3$ in case (P)---the case (P), $d=2$ is similar up to the
integrations by parts that are to be done with care) and we may define 
$L(t)=B(t)\,E(t)$ and, as for the Vlasov-Poisson system,
\beas
\frac{dL}{dt}
&=&
\alpha\ixed{V(t,x)\biggl(\dot{B}F
-\frac{\dot{R}}{R}B\partial_x(xF)\biggr)}\nonumber\\
&&
{}+\biggl(\dot{B} -
2\frac{\dot{R}}{R}B\biggr)\ixed{F(t,x,\eta)\vert\eta\vert^2}\label{LyapW}\\
&&
{}+\biggl[\dot{B}\frac{\ddot{R}}{R}+\biggl(\frac
d{dt}(\frac{\ddot{R}}{R})+2\frac{\ddot{R}}{R}\frac{\dot{R}}{R}\biggr)B\biggr]
\ixed{F(t,x,\eta)\vert x\vert^2}\; .
\eeas
 
In the case of the coupling with the Poisson equation ($d\geq 2$) , the
conditions on
$L$ that are sufficient for it to be nonincreasing are exactly the same as for
(R$'$VP) in the Poisson case (P): see Section 4. The detailed justifications of
the computations for initial data
\[ \vert \psi(t=0,\cdot)\vert^2=\int_{\R^d}w(t=0,\cdot,v)\; dv\in L^1(\R^d) \]
are not given here, and we shall refer to \cite{C1} for a proof if $d\geq 3$ 
in the context of the Schr\"odinger-Poisson system.

\begin{thm}\label{WP} 
Assume that $w$ is a solution of (WP) with $M=\Vert
w(t,\cdot,\cdot)\Vert$ and that $t\mapsto R(t)$ is the solution of Equation (\ref{req}): 
$\ddot R + \varepsilon R^{1-d}=0$, $R(0)=1$, $\dot{R}(0)=0$. The function $t\mapsto
L(t)$ defined above for $d\geq 2$ (with $B=R^{d-2}$) is decreasing for $d=2, 3$,
constant for $d=4$, and satisfies for any $d\geq 2$
\[ 
\frac {dL}{dt}= (d-4)\, \dot{R}R^{d-3}
\int\!\!\!\!\!\int_{\R^d\times\R^d} \biggl| v - \frac{\dot R}{R} x \biggr|^2
w\, dv\, dx\; . 
\]
In the plasma physics case $\varepsilon=-1$, $L$ is bounded from below and for
$d=2,3$,
\[ 
\int_0^{+\infty} \dot R(s)R^{d-3}(s)
\biggl(\int\!\!\!\!\!\int_{\R^d\times\R^d} 
\biggl| v - \frac{\dot R(s)}{R(s)} x \biggr|^2
w(s,x,v)\, dv\, dx\biggr)\, ds<+\infty\; . 
\]
\end{thm}
However, the results on the dispersion for the Vlasov-Poisson system cannot be
transposed straightforwardly because of the lack of positivity of $w$ and one
has to be very careful to recover the estimates given in \cite{IZL} for $d=3$.
In dimension $d=2$, the situation is even worse because the boundedness of $L$
from below is not obvious at all. In that sense, the Schr\"odinger formulation
of the problem is more suitable.

\subsection{Schr\"odinger equation}
The Lyapunov function for the Schr\"odinger equation is easily found 
by simply considering the Wigner transform.  However,
it is interesting to 
realize how the method of Section 4 applies directly. According to the Weyl
quantification and the Wigner transform, the
operator $i\hbar\partial_x$ corresponds to the variable $v$: 
the change of variables $\eta=v-\frac{\dot{R}}Rx$
therefore means that instead of $i\hbar\partial_x$ we consider the new
operator $i\hbar\partial_x-\frac{\dot{R}}Rx$:
\[ 
\phi\mapsto (i\hbar\partial_x-\frac{\dot{R}}Rx)\phi =e^{-i\frac{\dot{R}}R
\frac{\vert x\vert^2}{2\hbar}} i\hbar\partial_x
\biggl( e^{i\frac{\dot{R}}R\frac{\vert x\vert^2}{2\hbar}}\phi\biggr)\; .
\]
For that purpose, we may consider the new wave function $\phi(t,x)$ such that
\[ \psi(t,x)=e^{i\frac{\dot{R}}R\frac{\vert x\vert^2}{2\hbar}}\phi(t,x)
\; . \]
$\phi(t,x)$ solves the rescaled Schr\"odinger equation (R$'$S)
\[ i\hbar\partial_t\phi=-\frac 12 {\hbar}^2\Delta \phi+ 
(V+\frac{\ddot{R}}{2R}|x|^2)\phi
-\frac{i\hbar\dot{R}}{2R}(d\phi+2x\cdot\partial_x\phi)\; .\] 
If we define the
potential energy term by $W[\phi]=2V\vert\phi\vert^2$, $W[\phi]= 
V\vert\phi\vert^2$, or $W[\phi]=\frac 2{p+1}\vert\phi\vert^{p+1}$ in case (L),
(P), or (NL) respectively, the corresponding energy 
is given by 
\[ 
E(t)=\ixd{\biggl(\hbar^2\vert\nabla\phi\vert^2+W[\phi]+\frac{\ddot{R}}R\vert
x\vert^2 |\phi|^2\biggr)} 
\]
if $d\geq 3$, and
\[ 
E(t)=\ixd{\biggl(\hbar^2\vert\nabla\phi\vert^2+W[\phi]+\frac{\ddot{R}}R\vert
x\vert^2|\phi|^2\biggr)} 
+
\log R(t)\frac 1{2\pi}\biggl(\ixd{\vert\phi\vert^2}\biggr)^2 
\]
if $d=2$ in case (P). We may then build the Lyapunov functional in
the same way as for the solution of the Wigner equation. Going back to the
original variables, we have to replace $|\nabla\phi |^2$ by $|(\nabla-i\frac
{\dot{R}}{\hbar R}x)\psi |^2$.

\medskip
The Schr\"odinger-Poisson system and its asymptotics has been studied in
\cite{DF1,DF2,IZL}. More recently, a theory for $L^2$ solutions corresponding
to mixed quantum states has been established by F.~Castella (see \cite
{C1,C2}). In the case of a pure quantum state, J.~L.~Lopez and J.~Soler in
\cite{LS1,LS2} also gave detailed results on the asymptotic behaviour using a
linear scaling approach in the continuation of the method developed S. Kamin 
and J.~L.~V\'azquez. The main interest of our approach is that it gives a refined
estimate for $d=3$ and is adapted to the limit case $d=2$ as well.

Concerning the notion of solution we may assume that it is as smooth as
desired and refer to \cite{C1,C2} for minimal requirements (estimates for weak
solutions are built using approximating smooth solutions).

\begin{thm}\label{SP} 
Assume that $d\geq 2$ and consider a solution of the
Schr\"odinger-Poisson system. With the above notation
\[ 
L(t)= R^{d-2}(t)\int
\biggl|(\nabla-i\frac {\dot{R}}{\hbar R}x)\psi \biggr|^2\; dx +
R^{d-2}(t)\int V\; |\psi|^2\; dx +\frac 1{R^2(t)}\int |x|^2\; |\psi|^2\; dx \]
for $d=3,4$, and
\[ 
L(t)=\int\biggl(
\biggl|(\nabla-i\frac {\dot{R}}{\hbar R}x)\psi \biggr|^2+ V\;
|\psi|^2 +\frac {|x|^2}{R^2(t)}\; |\psi|^2\biggr)\; dx +
\frac{\Vert\psi(t,\cdot)\Vert^2_{L^2}}{2\pi}\log R(t)\]
for $d=2$, is decreasing for $d=2$, $3$ and constant for $d=4$ if $t\mapsto
R(t)$ is a solution of $\ddot R =  R^{1-d}$, $R(0)=1$, $\dot{R}(0)=0$.
As a consequence $n(t,x)= |\psi(t,x)|^2$ is decreasing: there exists a
constant $C>0$ such that 
\be 
\Vert n(t,\cdot)\Vert_{L^p(\R^d)}\leq C\cdot {\dot{R}}^{d(\frac 2p
-1)}\cdot R^{d(\frac 1p-\frac 12)(\frac d2-1)} \label{interpSP}
\ee
for any $p\in[2,\frac{2d}{d-2}]$ if $d=3$, $4$ and
$\liminf_{t\to +\infty}\Vert n(t,\cdot)\Vert_{L^p(\R^2)}=0 $ for 
$p\in ]2,+\infty[$ if $d=2$.
\end{thm}

Note that for $d=3$, $p=10/3$, we recover the same exponents as for the
Vlasov-Poisson system. For $d=2$, exactly the same estimate as in the proof of
Theorem \ref{LyaVPp2-4} holds:
\[ 
\liminf_{t\to +\infty}\int_{\R^2}
\biggl|(\nabla-i\frac {\dot{R}}{\hbar R}x)\psi \biggr|^2\; dx=0\; . 
\]
The crucial ingredient in the proof of this theorem is the following interpolation lemma 
(see \cite{DF1,DF2} and \cite[Cor.~5.5]{IZL}) 
which plays a role similar
to the one of Equation \eqn{interpolation} for the Vlasov-Poisson system:

\begin{lemma}\label{ILZlemma} Assume that $d\geq 3$. There exists a constant
$C>0$ depending only on $d$ such that, for any $u\in H^1(\R^d)$ such that 
$x\mapsto xu(x)$ belongs to $L^2(\R^d)$,
\[ \Vert u\Vert_{L^p(\R^d)}\leq C\Vert u\Vert_{L^2(\R^d)}^a
\Vert (x+it\nabla)u\Vert_{L^2(\R^d)}^{1-a}\cdot t^{-(1-a)} \]
for any $p\in [2, \frac {2d}{d-2}]$, $a=\frac d2 (\frac 2p-\frac {d-2}d)$.
\end{lemma}
The proof of Lemma \ref{ILZlemma} is easily established using the
Gidas-Nirenberg inequality
\[ \Vert u\Vert_{L^p(\R^d)}\leq [C(d)]^{1-a}\Vert u\Vert_{L^2(\R^d)}^a
\Vert \nabla u\Vert_{L^{\frac {2d}{d-2}}(\R^d)}^{1-a} \]
where $C(d)$ is the Sobolev constant corresponding to the injection of 
$H^1(\R^d)$ into $L^{\frac {2d}{d-2}}(\R^d)$, and the decomposition $u=\rho
e^{i\varphi}$ which holds at least for smooth enough functions (the
conclusion holds by a density argument). We may then write 
\[ \Vert (x+it\nabla)u\Vert_{L^2}^2 = t^2\int{\vert\nabla\rho\vert^2}\; dx +
\int{\vert x\rho+t\rho\nabla\varphi\vert^2}\; dx \geq t^2\Vert\; \nabla\vert
u\vert\;\Vert_{L^2}^2\; ,\]
which proves the interpolation results.
 
{\bf Proof of Theorem~\ref{SP}}.  
One has to replace $1/t$ by $\dot{R}/R$: for $d=3$ or $4$,
we may refer to \cite{IZL} for the proof of Equation \eqn{interpSP}, where it
is done in the case $R=t$. For $d=2$, the argument is similar, the main
step being the proof of the boundedness of $L$ which goes exactly as in
the Vlasov-Poisson case (see the proof of Theorem~\ref{LyaVPp2-4}).\prfe
\medskip

We conclude this section by considering the case of the Nonlinear
Schr\"odinger equation which allows us to make an explicit link with the
pseudo-conformal law. If $W[\phi]=\frac 2{p+1}\vert\phi\vert^{p+1}$, a direct
computation gives 
\[ 
\frac {dE}{dt}=-\frac {d(p-1)\dot{R}}{2R}\int\!\! W[\phi]\, dx
-2\frac{\dot{R}}R\hbar^2\int\!\!|\nabla\phi|^2 dx +
\biggl(\frac d{dt}(\frac{\ddot{R}}{R})
+2\frac{\ddot{R}}{R}\frac{\dot{R}}{R}\biggr)
\int\!\!{\vert x\vert^2 |\phi |^2} dx 
\]
and $L(t)=B(t)\,E(t)$ is decreasing if $B(t)=R^q(t)$, $q=\min((p-1)d/2, 2)$,
$\ddot{R}=1/R^{q+1}$, $R(0)=1$, $\dot{R}(0)=0$. In the next result we
are again
not interested in the weakest notion of solution and assume that the
solution is global in $t$ and as smooth and sufficiently decreasing at
spatial infinity as necessary to justify any integration by parts 
in the computations.

\begin{thm}\label{NLS} Assume that $d\geq 2$ and consider a global solution of
the Nonlinear Schr\"odinger equation (NLS)
\[
i\hbar\partial_t\psi=-\frac{\hbar^2}2\Delta\psi +|\psi
|^{p-1}\psi\; .
\] 
Then with the above notation
\be L(t)= R^q(t)\int\biggl(|(\nabla-i\frac {\dot{R}}{\hbar R}x)\psi |^2+ \frac
2{p+1} |\psi|^{p+1}\biggr)\; dx + \frac 1{R^2(t)}\int |x|^2\;
|\psi|^2\; dx \label{PsCfLaw}\ee
is decreasing. \end{thm}

Decay estimates can of course be deduced from Lemma \ref{ILZlemma} as for the
Schr\"odinger-Poisson system. The details of the computations for the proof of
Theorem \ref{NLS} are left to the reader.
\smallskip

A simple method to understand the pseudo-conformal 
law is simply to look for a pseudo-conformal invariance of the equation, 
{\sl i.e.\/} a transformation which leaves the equation invariant. If 
$u(t,x)$ is a solution of (NLS) in the focusing or in the defocusing case
($\varepsilon=-1$):
\[ i\hbar\partial_t u=-\frac 12 {\hbar}^2\Delta u+ \varepsilon
\vert u\vert^{p-1} u\; , \]
let us look for a function $t\mapsto (R(t), \tau(t), \omega(t))$ such that 
$(\tau, \xi)\mapsto v(\tau, \xi)$ given by 
\[ u(t,x)=\frac 1{R^\alpha(t)}e^{i\omega(t)\frac {\vert x\vert^2}2}
v(\tau(t), \xi(t))\;, \quad\xi(t)=\frac x{R(t)} \]
is a solution of (NLS) too for some $\alpha\in\R$. It turns out that this 
is possible only in the case $p-1=\frac 4d$ (critical case), and in 
that case, $t\mapsto (R(t), \tau(t), \omega(t))$ solves the system
\[ \frac {d\tau}{dt}
=\frac 1{R^2}\; ,\quad\frac {dR}{dt}=2\omega R\; ,\quad \frac
{d\omega}{dt}=-2\omega^2\; . \]
The solution is given by
\[ \omega(t)=\frac{\omega_0}{1+2\omega_0t}\; ,\quad R(t)= R_0(1+2\omega_0t)\;
,\quad\tau(t)=\frac t{R_0^2(1+2\omega_0t)}+\tau_0\; .\]
This transformation can be found in \cite{Me} (see also \cite{KW} for 
instance). The conservation of the energy after rescaling (conservation of the
energy for $v$) gives the following conservation law for $u$:
\be \frac d{dt}\biggl(R^2(t)\int_{\R^d}\vert\nabla u(t,x)-i\omega(t)x
u(t,x)\vert^2 \; dx -\frac{d\varepsilon}{d+2}\int_{\R^d}\vert u(t,x)\vert^{\frac
2d(d+2)} \; dx\biggr)=0\; . \label{PCL}\ee
This expression clearly corresponds to the case
$q=2=(p-1)d/2$, and the pseudo-conformal law is nothing else than the
expression of $dL/dt$ where $L$ is given by Equation \eqn{PsCfLaw}. As we
already noticed already several times, one may replace $\omega(t)$ and $R(t)$
by their equivalents as $t\to +\infty$, which is the same as considering the
singular solution corresponding to the limit $\omega_0\to +\infty$ and
$R_0\omega_0\to 1$, and recover instead of Equation \eqn{PCL} the more
classical form for the conformal invariance law:
\[ 
\frac d{dt}\biggl(t^2\int_{\R^d}
\biggl|\nabla u(t,x)-i\frac x{2t} u(t,x)\biggr|^2
\; dx -
\frac{d\varepsilon}{d+2}\int_{\R^d} u(t,x)\vert^{\frac 2d(d+2)}
\; dx\biggr)=0\; . \]

\medskip\noindent{\bf Acknowledgements.} The first author wants to thank 
Fran\c{c}ois Castella who introduced him to the question of the 
dispersion in the Schr\"odinger-Poisson system, and pointed out that the method 
should adapt to the case $d=2$ and Jalal Shatah for stimulating discussions
on the Nonlinear Schr\"odinger equation. He also thanks the Mathematisches
Institut der Universit\"at M\"unchen for welcoming him, and the TMR network No.
ERB FMRX CT97 0157 and the Erwin Schr\"odinger Institute of Wien for partial
support.


\end{document}